\documentclass[psamsfonts]{amsart}

\newtheorem{theorem}{Theorem}[]
\newtheorem{lemma}[theorem]{Lemma}
\newtheorem{proposition}[theorem]{Proposition}
\newtheorem{corollary}[theorem]{Corollary}

\theoremstyle{definition}

\newtheorem{example}[theorem]{Example}

\newtheorem{remarks}[theorem]{Remarks}

\theoremstyle{remark}

\numberwithin{equation}{section}

\newcommand{\smano} {\smallskip \noindent}

\newcommand{\ZZ}{\mathbb{Z}}

\newcommand{\CC}{\mathbb{C}}

\newcommand{\PP}{\mathbb{P}}

\newcommand  {\lra}     {\longrightarrow}

\renewcommand{\O}       {\mathcal{O}}

\newcommand  {\ra}      {\rightarrow}

\def\mydate{\number\day\space\ifcase\month \or January\or February\or March\or April\or May\or
June\or July\or August\or September\or October\or November\or
December\fi \space\number\year}

\usepackage{amscd}
\usepackage{amssymb}

\begin{document}

\title[]{A class of
counter-examples to the hypersection problem based on forcing
equations}

\author[Holger Brenner]{Holger Brenner}
\address{Mathematische Fakult\"at, Ruhr-Universit\"at, Universit\"atsstr. 150,
               44780 Bochum, Germany}

\email{brenner@cobra.ruhr-uni-bochum.de}

\begin{abstract}
We give a class of three-dimensional Stein spaces $W$ together
with a hypersurface $H$, such that the complement $W-H$ is not
Stein, but such that for every analytic surface $S \subset W$
the complement $S-S \cap H$ is Stein. This class is constructed using
forcing equations and gives new counter-examples to the hypersection
problem.
\end{abstract}

\maketitle


%===========================================================

\bigskip
\noindent
{\bf 0. Introduction.}
Let $W$ be a complex
Stein space of dimension $\geq 3$
and let $H \subset W$ be an analytic hypersurface, $U=W-H$.
Suppose that for every analytic hypersurface $S \subset W$
the intersection $U \cap S$ is Stein,  is then $U$ itself Stein?
This question is called the hypersection problem, see
\cite{diederich} for a general treatment and related problems.
The first counter-example to this question was given
by Coltoiu and Diederich in \cite{coltoiu1}, using the affine cone
over the complement of two sections
on some ruled surface over an elliptic curve.
In this way they get a normal three-dimensional isolated
singularity.

In this paper we present another class of three-dimensional
Stein spaces $W$ together with a hypersurface $H$
fulfilling the assumptions in the hypersection problem, but not
its conclusion.
The class is constructed in the following way:
we start with a two-dimensional normal
affine cone $X$ over a smooth projective
curve and the vertex point $P \in X$. Suppose that we have three homogeneous
functions $f_1,f_2$ and $f_0$ on $X$.
Then, under suitable conditions,
$W=V(f_1t_1+f_2t_2+f_0) \subset X \times \CC^2$ and
the hypersurface $H=p^{-1}(P)$ have the desired properties, see
Theorem \ref{theorem}.
These conditions reduce to numerical conditions,
which are easily to verify, see Corollary \ref{corex}
and Example \ref{example}.

\bigskip
\noindent
{\bf 1. Forcing equations.}
Let $X$ be an irreducible normal complex Stein space of dimension $d$
together with a point
$P \in X$.
Let $f_1, \ldots ,f_n$ be holomorphic functions on $X$
such that the common zero set of these functions is exactly $P$.
Let $f_0$ be another holomorphic function vanishing at $P$.
Then we consider the complex space $W \subset X \times \CC^n$ defined
by the equation $f_1t_1+ \ldots +f_nt_n +f_0 =0$,
$$W 
=\{ (x,t_1, \ldots ,t_n) \in X \times \CC^n:
\, f_1(x)t_1+ \ldots + f_n(x)t_n + f_0(x)=0 \} \, .$$
The equation $f_1t_1+ \ldots +f_nt_n +f_0=0$ is called a forcing
equation, since it forces $f_0$ to lie in the ideal
generated by $f_1, \ldots ,f_n$.
Forcing equations and the algebras defined by them play an important
role in the theory of closure operations for ideals, e.g.
tight closure and solid closure, see \cite{hochster1}.
Let $p: W \ra X$ be the projection.
If $d \geq 2$, then $n \geq 2$ and
$W$ is an irreducible Stein space of dimension $d+n-1$.
Let $U_i= \{ x \in X:\, f_i(x) \neq 0 \}$ and $U= \bigcup_{i=1}^n U_i =X-P$.
Resolving $t_i$ shows that $p^{-1}(U_i) \cong U_i \times \CC^{n-1}$.
The transition mappings however are only affine-linear, not linear,
hence $W|_U$ is not a vector bundle.
$H:= p^{-1}(P) \cong \CC^n $ is a closed subset in $W$, which is a
hypersurface in case $d=2$.
The existence of a section $X \ra W$ is equivalent with
$f_0 \in (f_1, \ldots ,f_n)$ over $X$.

\begin{lemma}
\label{superhoehe}
Let $X$ be a normal irreducible
Stein space together with a point $P$ and let
$f_1,\ldots, f_n,f_0 \in \Gamma(X,\O_X)$ be holomorphic functions
on $X$.
Let $X'$ be another irreducible Stein space of the same
dimension and
let
$\psi: X' \ra X$ be a holomorphic mapping
such that $\psi^{-1}(P)$ contains isolated points.
Suppose that
$f_0 \circ \psi \in (f_1 \circ \psi, \ldots , f_n \circ \psi)$ in
$\Gamma(X',\O_{X'})$.
Then
$f_0 \in (f_1, \ldots, f_n)$ in $\O_{X,P}$.
\end{lemma}
\proof
Let $Q$ be an isolated point over $P$.
Then there exist open neighborhoods $Q \in U$ and $P \in V$
such that $\psi : U \ra V$ is finite, see \cite{CAS}, Ch. 3.2.
Due to the finite mapping theorem, $\psi _* (\O_{X'})$ is a coherent
analytic algebra on $P \in V \subseteq X$.
Furthermore it is torsionfree due to the assumptions on the dimension.
Since $X$ is normal,
we have the trace map $tr: \psi _* (\O_{U}) \ra \O_V$,
which gives the result.
\qed

\begin{corollary}
\label{superhoehe2}
Let $X$ be a normal irreducible
Stein space of dimension $d$ together with a point $P$ and let
$f_1,\ldots, f_n,f_0 \in \Gamma(X,\O_X)$ be holomorphic functions such that
$f_0 \not \in (f_1, \ldots ,f_n) \O_P$. Let
$$W= V(f_1t_1+ \ldots +f_nt_n +f_0) \subset X \times \CC^n
\stackrel{p}{\lra} X\, .$$
and $H=p^{-1}(P)$.
Let $T$ be an irreducible complex space of dimension $d$
and let $\varphi: T \ra W$ be a holomorphic map.
Then $\varphi^{-1}(H) \subset T$ contains no isolated points
and the codimension of $\varphi^{-1}(H)$
is $\leq d-1$.
\end{corollary}

\proof
We look at the composed mapping
$\psi =p \circ \varphi :T \stackrel{\varphi} {\ra} W \stackrel{p}{\ra} X$.
Since it factors through $W$ it follows that
$f_0 \circ \psi \in (f_1 \circ \psi, \ldots , f_n \circ \psi)$
in $\Gamma(T,\O_T)$. Due to the Lemma
$\psi^{-1} (P)= \varphi^{-1}(H)$ cannot contain
isolated points.
\qed

\medskip
\noindent
With this result we can establish the hypothesis of the hypersection problem
in a broad class of example where the base space $X$
is two-dimensional, $n=2$ and $H=p^{-1}(P) \cong \CC^2$
is a hypersurface in three-dimensional $W$.

\begin{proposition}
\label{hypersechypo}
Let $X$ be a normal irreducible two-dimensional
Stein space together with a point $P$.
Let $f_1,f_2, f_0 \in \Gamma(X,\O_X)$ be holomorphic functions such that
$P= \{ f_1=f_2= 0\}$ and suppose that
$f_0 \not\in (f_1, f_2) \O_P$. Let
$$W= Z(f_1t_1+f_2t_2+ f_0) \subset X \times \CC^2 \stackrel{p}{\lra} X \, .$$
and let $H=p^{-1}(P) \subset W$.
Then for every analytic surface $S \subset W$ the complement of
$S \cap H \subset S$ is Stein.
\end{proposition}
\proof
We may assume that $S$ is irreducible, let $\tilde{S}$ be its
normalization and let
$\varphi : \tilde{S} \ra W$ be the corresponding mapping.
Due to Cor. \ref{superhoehe2}
we know that
$\varphi^{-1}(H) \subset \tilde{S}$ contains no isolated points.
Hence
$\varphi^{-1}(H)$ is a pure curve on a normal Stein surface
and
due to \cite{simha} its complement is Stein.
But then also $S -S \cap H$ itself is Stein.
\qed

\bigskip\noindent
{\bf 2. The graded situation.}
We have to look for examples of the type described in Proposition
\ref{hypersechypo}
where $W-H$ is not Stein.
To this end we look at the graded situation.
Let $Y \subseteq \PP^N$ be a smooth projective variety
with the very ample line bundle $H_Y \ra Y$
(which is the restriction of $\O_{\PP^N}(1)$ to $Y$)
and let $X \subseteq \CC^{N+1}$ be the corresponding
affine cone.
Let $P$ be the vertex of the cone and assume that $X$ is normal.
Recall that we have an action of $\CC^*$ on $X$, which is free on $U=X-P$.
$Y$ is the quotient of this action and $U \ra Y$ is a $\CC^*$-principal
bundle.
A number $e \in \ZZ$ defines the action on $X \times \CC$ by
$\lambda (x,t):=(\lambda x, \lambda ^e t)$,
this action is free over $U$ and the quotient is
the line bundle $H_Y^e \ra Y$.

Suppose that the holomorphic functions
$f_i$ are homogeneous of degree $d_i$, i.e.
$f_i(\lambda x )= \lambda^{d_i}f(x),\, x \in X,\, \lambda \in \CC^*$.
We may consider
a homogeneous holomorphic function $f$ of degree $d$
as a section $Y \ra H_Y^d$ and
as a mapping of line bundles $H_Y^{-d} \ra Y \times \CC$ or
$H_Y^{e} \ra H_Y^{e+d}$.

\begin{proposition}
\label{buendelalsproj}
Let $X$ be a normal affine cone over a smooth
projective variety $Y$ and let $P$ be the vertex point.
Let $f_1,\ldots,f_n$ be homogeneous functions
such that $P= \{ f_1= \ldots =f_n=0 \}$, $U=X-P$.
Let $d_i$ be the degrees of $f_i$ and let $e_i$ numbers such that
$m=d_i+e_i$ is constant.
Then the following hold.

\renewcommand{\labelenumi}{(\roman{enumi})}
\begin{enumerate}

\item
There is an exact sequence of vector bundles
$$0 \longrightarrow V_m \longrightarrow
H_Y^{ e_1} \times_Y \ldots \times_Y H_Y^{ e_n}
\stackrel{\sum f_i}{\lra} H_Y^{ m } \longrightarrow 0 \, .$$

\item
$V_m$ is the quotient of $V(f_1t_1+ \ldots + f_nt_n)|_U$
by the action of $\CC^*$
given by
$\lambda (x,t_1, \ldots, t_n)
=(\lambda x, \lambda^{e_1}t_1, \ldots ,\lambda^{e_n}t_n)$.

\item
We have ${\rm Det}\, V_m \cong H_Y^k$, where
$k= \sum_{i=1}^n e_i -m= -\sum_{i=1}^n d_i +(n-1)m $.

\item
$V_{m'}= V_m \otimes H_Y^{(m'-m)}$.
The projective bundle
$\PP(V_m)$ is independent of the choosen degree $m$.

\end{enumerate}
\end{proposition}

\proof
(i).
We consider the $f_i$ as morphisms
$H_Y^{e_i} \ra H_Y^m$. The morphism of vector bundles
$H_Y^{ e_1} \times_Y \ldots \times_Y H_Y^{ e_n}
\stackrel{\sum f_i}{\ra} H_Y^{ m } $
over $Y$ is surjective, because the $f_i$ do not have a common zero on $Y$.
Hence the kernel is a vector bundle $V_m$ on $Y$ of rank $n-1$.

\smano
(ii). The pull back under $q: U \ra Y$ of the exact
sequence in (i) gives
$$
0 \longrightarrow q^* V_m \longrightarrow
U \times \CC^n
\stackrel{\sum f_i}{\lra} U \times \CC \longrightarrow 0 \,$$
together with the described action,
and $q^*V_m= V(f_1t_1+ \ldots +f_nt_n)|_U$.
(iii) and (iv) follow.
\qed

\begin{proposition}
\label{buendelsequenz}
Let $X$ be a normal affine cone over a smooth
projective variety $Y$ and let $P$ be the vertex point, $U=X-P$.
Let $f_1,\ldots,f_n$ be homogeneous functions
such that $P= \{ f_1= \ldots =f_n=0 \}$.
Let $f_0$ be another homogeneous function, $d_i = \deg (f_i)$,
and let $e_i$ be numbers such that
$m=d_i+e_i$ is constant for $i=0, \ldots ,n$.
Let $V_m$ {\rm (}resp. $V'_m${\rm )}
be the vector bundle on $Y$
defined in Proposition {\rm \ref{buendelalsproj}} with respect to 
$f_1, \ldots, f_n$
{\rm (}resp. $f_0, \ldots, f_n${\rm )}.
Then the following hold.
\renewcommand{\labelenumi}{(\roman{enumi})}
\begin{enumerate}

\item 
There is an exact sequence of vector bundles on $Y$:
$$0 \longrightarrow V_m \longrightarrow V'_m
\stackrel{t_0}{\longrightarrow} H_Y^{ e_0} \longrightarrow 0 \, .$$

\item
The corresponding embedding $\PP(V_m) \hookrightarrow \PP(V'_m)$ is
independent of $m$, $\PP(V_m)$ is a divisor on $\PP(V'_m)$.

\item
Let $e_0=0$. The normal bundle for $\PP(V) \hookrightarrow \PP(V')$
on $\PP(V)$
is $H_{\PP(V)}$, where $H_{\PP(V)}$ denotes the relative very ample
line bundle on $\PP(V)$.

\item
Let $e_0=0$.
$W|_U \longrightarrow \PP(V') - \PP(V) $ is a quotient
of the action on $W=V(f_1t_1+ \ldots +f_nt_n+f_0)$
given by
$$
\lambda (x,t_1, \ldots, t_n)
\mapsto (\lambda x, \lambda ^{e_1}t_1,\ldots , \lambda^{e_n}t_n) \, .
$$

\end{enumerate}
\end{proposition}
\proof
(i).
The mappings in the sequence follow from
the defining sequences for $V_m$ and $V'_m$.
The exactness of the sequence follows from diagramm chasing.
(ii) is clear.

\smano
(iii).
Since we assume $e_0=0$, $t_0$ is a global function on $V'$ and it
is a global section in the relative very ample line bundle $H_{\PP(V')}$
on $\PP(V')$, and $\PP(V)$ is the corresponding divisor.
Therefore the normal bundle of this embedding
is $i^*H_{\PP(V')} \cong H_{\PP(V)}$.

\smano
(iv).
First we may identify the closed subset
$\{ Q \in V':\, t_0(Q) =1 \}$ with $\PP(V') -\PP(V)$.
The described action on
$W $
respects the forcing equation, for
$ f_1(\lambda x) \lambda^{e_1} t_1 +
\ldots + f_n(\lambda x) \lambda^{e_n} t_n +f_0(\lambda x)
= \lambda^{d_1} \lambda^{e_1} f_1(x)t_1+ \ldots +
\lambda^{d_n} \lambda^{e_n} f_n(x)t_n +\lambda^{d_0} f_0(x)
=\lambda^m (f_1(x)t_1+ \ldots + f_n(x)t_n +f_0(x)) =0$.
This action on $W|_U =W-p^{-1}(P)$
is the same action as the action
on the vector bundle
$V(f_1t_1+ \ldots +f_0t_0)|_U$ restricted to $t_0=1$ described in
Proposition \ref{buendelalsproj}(ii).
Its quotient is $\{ Q \in V': \, t_0(Q)=1\}$.
\qed

\medskip\noindent
Now we specialize to the two-dimensional situation.

\begin{corollary}
\label{intersection}
Let $X$ be a normal affine two-dimensional cone over a smooth
projective curve $Y$ and let $P$ be the vertex point.
Let $f_1,f_2$ be homogeneous functions
such that $P= \{ f_1=f_2=0 \}$.
Let $f_0$ be another homogeneous function, $d_i = \deg (f_i)$,
and let $e_i$ numbers such that
$m=d_i+e_i$ is constant for $i=0,1,2$.
Let $V_m$ {\rm (}$V'_m${\rm )} be the corresponding vector bundles on $Y$.
Then the following hold.

\renewcommand{\labelenumi}{(\roman{enumi})}
\begin{enumerate}

\item
$\PP(V')$ is a ruled surface and $\PP(V) \subset \PP(V')$
is a section {\rm (}independent of $m${\rm )}.

\item 
We have
$V_m \cong H_Y^{e_0+d_0-d_1-d_2}$ and
the exact sequence
$$ 0 \longrightarrow H_Y^{e_0 +d_0-d_1-d_2} \longrightarrow V'
\longrightarrow H_Y^{ e_0} \longrightarrow 0 \, .$$

\item
Let $e_0=0$.
The normal bundel of the embedding
$Y \cong \PP(V) \subset \PP(V')$
is $H_Y^{d_1+d_2-d_0}$.

\item
The self intersection number of
$Y \cong \PP(V) \hookrightarrow \PP(V')$
is
$(d_1+d_2-d_0) \deg H_Y$.
\end{enumerate}
\end{corollary}

\proof
(i) is clear due to Proposition \ref{buendelsequenz}.

\smano
(ii)
From the defining sequence in Proposition \ref{buendelalsproj}
it follows that we have
$V_m \cong H_Y^{e_1} \otimes H_Y^{e_2} \otimes H_Y^{-m}
=H_Y^{e_1+e_2-m}=H_Y^{e_0+d_0-d_1-d_2}$.

\smano
(iii).
The normal bundel on $\PP(V)$ is $H_{\PP(V)}$
due to Proposition \ref{buendelsequenz}.
But for a line bundle this is just the negative
tautological bundel $-V$, therefore
$N=-V=H^{d_1+d_2-d_0}$.

\smano
(iv).
The self intersection number
is $\PP(V)^2= \deg_Y N= \deg_Y H^{d_1+d_2-d_0}=(d_1+d_2-d_0) \deg H_Y$.
\qed

\bigskip\noindent
{\bf 3. A class of examples.}

\begin{theorem}
\label{theorem}
Let $X$ be a normal affine two-dimensional
cone with vertex point $P$
over a smooth projective curve $Y$,
let $f_1,f_2$ and $f_0$ be homogeneous holomorphic functions on $X$
with degrees $d_1,d_2,d_0$
such that
$$(1)\,\, \, V(f_1,f_2)=P,\, \, \, \, \,
(2)\, \, \, f_0 \not\in (f_1,f_2) \O_{X,P}\, \, \mbox{ and }\, \, \, \,
(3) \, \, \, d_1 + d_2 - d_0 < 0 \, .$$
Let $W=V(f_1t_1+f_2t_2+f_0) \subset X \times \CC^2$ and $H=p^{-1}(P) \subset W$.
Then
$W-H \subset W  $
is not Stein,
but for every analytic surface $S \subset W$
the intersection $(W-H) \cap S$ is Stein.
\end{theorem}

\proof
We have to show that $W-H$ is not Stein.
Since $W-H=W|_U$, $U=X-P$, the quotient of this open subset under
the action of $\CC^*$ is $\PP(V')-\PP(V)$.
Due to \cite{matsushima} it is enough to show that this
complement of the section in the ruled surface $\PP(V')$
is  not Stein. The self intersection of $\PP(V)$ is
$( d_0 + d_1 -d_0) \deg H_Y < 0$,
hence $\PP(V)$ is due to
\cite{grauertmod}
contractible and the complement cannot be Stein.
\qed

\begin{corollary}
\label{corex}
Let $h$ be an irreducible homogeneous polynomial of degree $r $
in the three variables $x,y,z$ and suppose that $x,y$
are homogeneous parameters
{\rm (}i. e. that $h=c z^r +$ other terms, $c \neq 0${\rm )}
and let $X=V(h) \subset \CC^3$.
Suppose that $X$ is normal {\rm (}hence $Y$ is smooth{\rm )}.
Let $d_1, \,d_2 \geq 1$ and $d_0$ be degrees such that
$d_1 +d_2 <d_0 < r$.
Then $f_1=x^{d_1}, \, f_2 =y^{d_2}$ and $f_0 =z^{d_0}$
fulfill the conditions of the theorem.
\end{corollary}
\proof
We have to show that $z^{d_0} \not\in (x^{d_1},y^{d_2}) \O_P$.
For this we look at the completion
of the local ring,
$\hat{ O}_P /(x^{d_1},y^{d_2}) = \CC[[x,y,z]]/(h,x^{d_1},y^{d_2})$.
Since $d_0 < r= \deg (h)$ we see that $z^{d_0} \neq 0$ in this residue
class ring.
On the other hand,
the self intersection number is $(d_1+d_2-d_0) \deg \, H_Y < 0$.
\qed

\begin{example}
\label{example}
Let
$X=V(x^r+y^r+z^r) \subset \CC^3, \, r \geq 4$
be a Fermat type hypersurface, let
$$W= V(x^r+y^r+z^r,xt_1+yt_2+z^s) \subset \CC^5,
\, 3 \leq s < r \, \mbox{ and }\, H=V(x,y)\, .$$
Then the conditions in the corollary are fulfilled.

\smano
The easiest example of this type is the Fermat quartic $x^4+y^4+z^4$
together with $f_1=x,\, f_2=y$ and $f_0=z^3$.
Therefore
$$W= V(x^4+y^4+z^4,xt_1+yt_2+z^3) \, \mbox{ and }\, H=V(x,y) $$
gives an counter-example to the hypersection problem.
\end{example}

\begin{remarks}
The hypersurface $H$ in our example is the singular locus of $W$.
Since the normalization does not change the complement of $H$ and since
its preimage is still a hypersurface due to Corollary \ref{superhoehe2},
we also may get normal examples.

The condition $f_0 \not\in (f_1,f_2) \O_P$ in Corollary \ref{corex}
ensures that the divisor $\PP(V) \subset \PP(V')$
intersects every curve $C \subset \PP(V')$
positively. For a disjoined curve would yield a
closed punctured surface (its cone)
inside $W-H$.
If additionally $d_1+d_2 -d_0 >0$, then $\PP(V)$ is
an ample divisor and its complement is affine,
hence Stein.
What happens if $d_1+d_2-d_0=0$?
Then the complement is not affine, but it may be Stein.
For $h=x^3+y^3+z^3=0,\, f_1=x, \, f_2=y$ and $f_0=z^2$
we get an instance of the classical construction of Serre of
a Stein, but non-affine variety,
see \cite{umemura}.

\end{remarks}

%===========================================================

\end{document}